\documentclass[12pt,a4paper,oneside]{amsart}
\usepackage[a4paper]{geometry}
\ifx\counterwithout\undefined
\usepackage{chngcntr}
\fi

\usepackage[perpage,symbol*]{footmisc}

\counterwithout{subsection}{section}

\usepackage{mathptmx} 
\DeclareMathAlphabet{\mathcal}{OMS}{cmsy}{m}{n} 

\input{cyracc.def}
\font\tencyr=wncyr10
\def\rus{\tencyr\cyracc}

\DefineFNsymbolsTM{yoyo}{
  \textdagger    \dagger
  \textdaggerdbl \ddagger
  \textsection   \mathsection
  \textparagraph \mathparagraph
  \textbardbl    \|%
  {\textasteriskcentered\textasteriskcentered}{**}%
  {\textdagger\textdagger}{\dagger\dagger}%
  {\textdaggerdbl\textdaggerdbl}{\ddagger\ddagger}%
}
\setfnsymbol{yoyo}

\begin{document}

\title{Mathematicians going east}
\author{Pasha Zusmanovich}
\address{University of Ostrava, Ostrava, Czech Republic}
\email{pasha.zusmanovich@osu.cz}
\date{First written May 1, 2018; last revised May 28, 2018}
\thanks{arXiv:1805.00242}

\begin{abstract}
We survey emigration of mathematicians from Europe, before and during WWII,
to Russia. The emigration started at the end of 1920s, the time of 
``Great Turn'', and accelerated in 1930s, after the introduction in Germany of 
the ``non-Aryan laws''. Not everyone who wanted to emigrate managed to do so, 
and most of those who did, spent a relatively short time in Russia, being 
murdered, deported, or fleeing the Russian regime. After 1937, the year of 
``Great Purge'', only handful of emigrant mathematicians remained, and even less
managed to leave a trace in the scientific milieu of their new country of 
residence. The last batch of emigrants came with the beginning of WWII, when 
people were fleeing eastwards the advancing German army.
\end{abstract}

\maketitle

\section*{Introduction}
A lot is written about emigration of scientists, and mathematicians in 
particular, from Germany and other European countries between the two world 
wars. First and foremost, one should mention a detailed study \cite{ss-fleeing},
concentrating mainly on emigration to US, but also briefly covering emigration 
to other countries; then there are a lot of papers concentrating on emigration 
to specific countries: \cite{bers}, \cite{us} (US), \cite{fletcher}, 
\cite{nossum-kotulek} (UK), \cite{rider} (both US and UK), \cite{denmark} 
(Denmark), and \cite{turk} (Turkey). In absolute figures, the largest number of 
mathematicians emigrated to US (over 100 by many accounts), while in relative 
numbers (say, proportional to the size of the population, or to the number of 
actively working mathematicians in the host country) the first place belongs to
UK.

Concerning emigration to Russia, besides separate accounts of the fate of 
individual mathematicians scattered over the literature (a recent interesting contribution being \cite{odinets}), there is only a very brief survey in \cite[pp.~132--135]{ss-fleeing}. The latter account contains a 
lot of interesting material which we frequently cite below, but it is by all 
means not complete, and, as far as we know, no detailed survey of emigration of 
mathematicians to Russia in the specified period exists. That is what we 
attempt to do here. We do not claim any originality, our account is 
prosopographic, and most of our sources are secondary and tertiary ones. A 
tremendous impact of European emigrants on American mathematics, and, to a 
lesser degree, on mathematics in other countries, is well known, and we find it 
interesting to compare the situation with those in Russia.

A few words what is understood under ``emigration of mathematicians before or 
during WWII'' are in order. As curious as it may sound, we need to clarify all 
the terms used: what is a mathematician, what is emigration, and what is 
``before and during WWII''. We treat the term ``mathematician'' liberally and 
inclusive; for example, we include in our considerations the logician 
\textbf{Chwistek}, the statistician \textbf{Gumbel}, the chess player 
\textbf{Lasker}, the physicist \textbf{Mathisson}, the electrical engineer 
\textbf{Pollaczek}, the mechanical engineer \textbf{Sadowsky}, and the 
philosopher \textbf{Wittgenstein}. While the main occupation of these men was 
outside mathematics proper, the mathematical component of their deeds, at least
in a certain period of their lives, was strong enough to consider them as 
mathematicians (even if by a purely formal criteria to have works listed in 
\emph{Jahrbuch \"uber Fortschritte der Mathematik} or \emph{Zentralblatt}). 
Still, we have to draw a line between mathematics and (theoretical) physics somewhere, and do not include in our considerations such 
theoretical physicists as Guido Beck, Max Born, Boris Podolsky, or Victor 
Weisskopf (but, for example, do include \textbf{Rosen}, a coauthor -- together 
with Einstein -- of Podolsky, or \textbf{Lustig}, an experimental physicist 
turned mathematician -- if a minor one -- \emph{after} the emigration)\footnote{
Emigration of physicists to Russia is briefly discussed in \cite[\S 2]{hoch}.
}. 
On the other 
hand, following the accepted by the community definition of a mathematician as a
someone who is doing (or did) mathematical research, we do not include persons 
with a mere mathematical education, but who did not publish a single paper, 
like, for example, a pole Henryk Gustaw Lauer, who, after emigration to Russia,
held a high post in Gosplan (the central Russian planning agency at the time), 
and was murdered in 1937.

In the turmoil European situation at the first half of XX century with 
drastically changing borders, there is no clear-cut definition of what should be
considered as emigration. We understand ``emigration'' in cultural terms, rather
than in political or bureaucratic ones. Some of our protagonists, such as 
\textbf{Grommer}, \textbf{Plessner}, and \textbf{Walfisz}, were born at the 
territories which at the time belonged to Russia, or even, like in the case of \textbf{Plessner}, at the time of settling in Russia had a Soviet citizenship. Still, as all their 
formative years, including secondary and higher education, were spent in the 
German-speaking and cultural environment, it is safe to count them as emigrants.
Also, we include the cases when, at the outbreak of the war between 
Russia and Germany, people were fleeing the advancing German army from the 
territories recently annexed by Russia (eastern Poland\footnote{
Or, western Ukraine, according to Ukrainian terminology.
}
and Baltic states). On the other hand, we do not consider the cases where a 
person just remained on these territories after annexation (like, for example, 
the Ukrainian mathematicians who remained in Lvov, and most of the 
mathematicians from the Baltic states and Bessarabia; however, the case of 
\textbf{Lewicki} emphasizes the complexity of the situation: he is listed below
under the cases of ``Emigrations that did not occur'', as he did not emigrate, 
despite not small efforts from the Ukrainian side, to the Russian part of 
Ukraine before WWII, but after the war he became a Soviet citizen by the fact of
annexation of eastern Poland). Neither do we consider the cases when a person 
was raised and educated in Russia, and then lived abroad, as, for example, the 
case of Alexander Chuprov, an eminent statistician, who in 1925 declined a call
for return to Russia. Victor Levin, who was born and raised in Russia, got the 
university education and the first academic post in Germany, was forced to leave
Germany in 1933, and after many adventures, returned to Russia in 1938, is a 
borderline case which we have decided not to include either. 

As in the case of emigration to the West, the most of emigrants came from 
Germany or German-occupied territories; and most of them were either Jewish, 
or of a left-wing political persuasion (or both). We have also included the 
cases when people emigrated as early as in the second half of 1920s, i.e. before
the rise of the Nazi regime in Germany: the emigrants were, again, either 
Jewish, for whom, due to official or unofficial antisemitic policies, it was impossible to find a 
suitable academic position in their home countries, or left-wing political 
activists. But we have to set a lower bound somewhere, so we exclude, for 
example, the appearance in Russia of Ernst Kolman\footnote{
Ernst Kolman (1892--1979): a native of Czech lands (then part of 
Austro-Hungary), a high-ranking Soviet bureaucrat, and a Marxist-Leninist 
``philosopher''. He actively participated virtually in all assaults of the 
Russian communist regime on mathematics in 1920--1930s, spied on behalf of 
Russian authorities on scientists during their trips abroad, and penned many 
secret denunciations to OGPU/NKVD (the Russian secret police) on various 
mathematicians. ``An authoritative person'' according to Kolmogorov, ``a dark 
angel of Soviet mathematics'' according to others. As controversial as it may 
sound, Kolman was an able mathematician. True, the most of his writings were a 
mix of mathematical gibberish and a blatant communist propaganda, and he was 
ridiculed at the 1932 International Mathematical Congress in Z\"urich after 
delivering a talk about ``Karl Marx's foundations of differential calculus'', but this does not invalidate the fact that he has a 
couple of genuine and interesting mathematical papers under his belt.
}\label{foot-kolman}
whose circumstances of ``emigration'' were entirely different anyway from
all the other cases considered here: he came to (pre-communist) Russia as a 
prisoner of war during WWI. 

On the other hand, we set the upper bound to (the end of) WWII (the 
chronologically last case of an emigration which did not occur considered by us
is that of \textbf{Banach} in 1945), not considering such interesting cases as 
(re)emigration of Lev Kaluzhnin in 1951, or forced labor of German applied 
mathematicians in the Russian atomic bomb project.

In the subsequent sections, broken down by year of emigration, we list all the 
known to us cases of emigration of mathematicians to Russia. After that, in the
section entitled {\sc Emigrations that did not occur}, we list all such 
``unsuccessful'' attempts. The latter list is very eclectic, ranging from 
unsuccessful many-years persistent efforts to emigrate to Russia, till a casual
mentioning by a third party about possibility of emigration of that or another 
person. In an attempt to get as complete picture as possible, we are trying to 
collect all, however small, such evidences, following a short {\sc Conclusion} 
at the end.

Naturally, while telling our stories, we encounter, besides the main 
protagonists (whose names are displayed in boldface), a lot of other people. We 
assume that the reader is reasonably versed in the history of mathematics and 
neighboring areas, such as theoretical physics, in Russia and elsewhere, so such
names as Pavel Aleksandrov, Einstein, Gelfand, Kolmogorov, Edmund Landau, Luzin,
Weyl, etc., do not require any explanation. On the other hand, the lesser 
figures, or persons outside of the mathematical world, usually warrant a brief 
footnote.

\section*{Year 1925}

\subsection{Stef\'ania Bauer (nee Szil\'ard)}
(Reference: \cite{bauer}, chiefly describing the fate of her husband, Ervin 
Bauer).

Born in 1898 in Gy\H{o}r (then Austro-Hungary), sister of 
\textbf{Karl Szil\'ard}. 

Graduated from the Budapest University in 1919, the same year married Ervin 
Bauer, a noted physician and biologist, and a staunch supporter of the 
short-lived Hungarian Soviet Republic. After the fall of the Republic in 1919, 
and ensuing ``white terror'', the couple fled the country and lived, under 
hardship conditions, in G\"ottingen, Prague, and Berlin. After invitation to 
Ervin Bauer by Semashko, the Russian Commissar (= minister) of Public Health of
that time, they settled in 1925 in Russia, ``the country of their dreams'', 
first in Moscow, and then in Leningrad.

Stef\'ania became a member of the communist party, and worked, together with her
husband, at the Leningrad branch of the Institute of Experimental Medicine. 
According to her contemporaries, she was an able mathematician, but her life 
apparently was dominated by the life and deeds of her husband, whom she helped 
in his work in theoretical biology. During her life, she published only two 
short papers -- one in Hungarian in 1917 in number theory, and another one in
1934 in \emph{Matematicheskii Sbornik} in complex analysis.

The couple was arrested in 1937, and murdered in 1938. Their two sons were 
separated, given different family names, and send to orphanages.

\section*{Year 1929}

\subsection{Celestyn Burstin}
(References: 
\cite[pp.~26--27]{belarus}, \cite[pp.~9--10 of the Russian edition]{volynskii},
\cite{moo}).

Born in 1888 in Tarnopol (Austro-Hungary, currently Ukraine), graduated from the
University of Vienna in 1911, and got a doctor degree in 1912. Being a Jew and
a communist, had difficulty to get an academic appointment in German-speaking 
lands, and moved to Minsk in 1929, where he was appointed as professor at the
Belarusian State University, and as director of the Institute of Mathematics of
the Belarusian National Academy of Sciences; shortly thereafter he was also
elected as a member of the Academy.

In 1925 he joined the Austrian communist party, and after emigration to Russia,
the Soviet communist party. In 1931, after a coup d'\'etat in the Moscow 
Mathematical Society, led by the communist zealots preaching the ``class 
character of science'', became a member of the Society's ruling board.

Arrested in 1937, died in prison in Minsk in 1938.

There are very few sources about Burstin. The references above contain his 
formal and short biography, as well as a very limited piece of information about
Burstin's administrative activity in the Belarusian academy.

Works in differential geometry, measure theory, and general algebraic systems.

\subsection{Mikolaj Czajkowski (Nikolai Andreevich Chaikovskii)}

(References: 
\cite{voznyak}, \cite{voznyak-voznyak}, \cite{voznyak-voznyak-0}, \cite{vozna}).

Born in 1887 in Berezhany (Austro-Hungary, currently Ukraine). Obtained, not
without difficulty, a doctor degree at Vienna in 1911 under Franz Mertens. 
Taught at high schools, and at the private Ukrainian underground university in 
Lvov in 1922--1924. Being an Ukrainian patriot, he suffered from real and 
perceived Polish oppression, and, while in Lvov, had no contacts whatsoever with
the flourishing at that time Lvov school of mathematics. As early as 1924, he 
attempted, via Mikhail Kravchuk\footnote{
Mikhail Kravchuk (1892--1942): a noted Russian-Ukrainian mathematician (of the 
Kravchuk polynomials), arrested by NKVD in 1938, died in 1942 in a 
concentration camp.
}, 
to get a post in the Russian part of Ukraine (``the Great Ukraine'', as he put 
it) with aspirations to build ``Ukrainian mathematics''. These efforts succeeded
only in 1929, when he got a position at one of institutes of higher learning in
Odessa. He started to work there enthusiastically; in letters to his father who
remained in Poland, he praised his life and working conditions, but at the same
time complained about ``Jewish dominance'' and his ``lack of knowledge of the 
Marxist-Leninist theory''. 

Czajkowski was arrested in 1933, and spent 10 years in Gulag (the Russian system
of concentration camps); meanwhile, his family was deported from Odessa. At 1956
he managed to return to Lvov, and later he was appointed as professor there. 
Died peacefully in 1970.

Czajkowski was a multifaceted person: he was a conductor in a professional 
choir, wrote popular science articles, as well as science fiction (apparently 
the first in this genre in the Ukrainian language).

He published a few papers on Galois theory and elementary number theory, but as
they appeared in Ukrainian and in obscure periodicals, these works seemingly 
remained unknown. Later in his life he has switched to history of mathematics, 
pedagogy, writing encyclopedia articles, developing Ukrainian mathematical 
terminology, translation of classical mathematical works to Ukrainian, etc.

\subsection{Felix Frankl}
(References: 
\cite[p.~325 of the English edition]{alexandrov-autobiog}, 
\cite[p.~120]{vucinich}, \cite[p.~64]{lapko-lyusternik}, \cite{moo}).

Born in Vienna in 1905, got a doctor degree in 1927 under Hans Hahn. Since 
1928, member of the Austrian Communist party. He befriended Pavel Aleksandrov 
during the latter's stay in Vienna, and asked him to facilitate his emigration 
to Russia. For that, Aleksandrov lobbied Otto Yulievich Schmidt\footnote{
Otto Yulievich Schmidt (1891--1956): a colorful, influential, and controversial
figure in the Russian science in 1920--1940s: mathematician (seminal works in 
group theory and cosmology), explorer of Arctic (Hero of the Soviet Union and a
honorary member of the New York-based Explorers Club), mountaineer, high-ranking
Soviet bureaucrat, chief of naval military expeditions, chief editor of 
\emph{Matematicheskii Sbornik} and a dozen of other periodicals, director of the
main state publishing house and a few academic institutes, professor at a number
of Moscow universities admired by students, expert in monetary and taxation 
policies, and lady's man. ``A Renaissance man'', according to Kolmogorov. Was 
fiercely attacked by Kolman and Yanovskaya (of whom see a footnote below) on political and 
ideological grounds.
}, 
with success\footnote{
According to Aleksandrov, Schmidt's first reaction was: ``We have enough 
communists of our own; let him stay in Vienna and start a revolution in 
Austria'', but then he relented and changed his mind.
}.

In 1929 Frankl emigrated to Russia, and immediately was commissioned with an 
important political task: to report on ``Soviet works on topology'', the first 
in the planned series of such reports on various branches of mathematics, in the
framework of the ongoing attempt to introduce a centrally coordinated, and 
subject to communist guidance, five-year planning of mathematics. Apparently, 
the report was too mathematically competent and void of communist phraseology to
please his political bosses; instead, it became a standard reference in Russian
topological books for years to come.

Frankl's first workplace in Russia was the so-called Communist Academy. At the 
first USSR mathematical symposium in 1930 he made a talk with a telling title 
``Dialectical logic and mathematics''. In 1931, as a result of coup d'\'etat 
in the Moscow Mathematical Society mentioned under \textbf{Burstin}, he, 
together with the latter, emerged as a member of the Society's ruling board. 
Since that time, and till 1935 he was also a member of the editorial board of 
\emph{Matematicheskii Sbornik}, a flagship Russian mathematical journal of that 
time.

Later Frankl has drastically changed his topic (topology), and worked in the 
area of mechanics, in various research institutes in Moscow\footnote{
It is reported in \cite[p.~56]{myshkis}, that in a certain American scientific
directory of that time, two F. Frankls were mentioned: one from Vienna working 
in topology, and another from Moscow working in aerodynamics.
}.

In 1950 he was expelled from the communist party (a grave punishment in Russia
of that period), and was deported to Frunze (today Bishkek) in Kyrgyzstan. He 
was unable to return to Moscow, and died in 1961 in Nalchik, a small city at the
outskirt of Russian empire. This last period of his life he worked in 
differential equations. During his forced work in province (Frunze and Nalchik)
he supervised many PhD theses.

\subsection{Jacob Grommer}
(References: 
\cite[pp.~64--66]{transcending}, \cite{einst-1917}, \cite{einst-ioffe}).

Born in 1879 in Brest-Litovsk (then Russia, currently Brest, Belarus).

At the young age, he suddenly became interested in mathematics, and in a short
time underwent a transformation from an uneducated Yiddish-speaking Jew from 
shtetl, ignorant in world science, culture, and anything else but Talmud, to
a brilliant doctoral student under Hilbert in G\"ottingen\footnote{
``If students without the gymnasium diploma will always write such dissertations
as Grommer's, it will be necessary to make a law forbidding the taking of the 
examination for the diploma'', reportedly said Hilbert 
(\cite[p.~143 of the 1996 edition]{reid-hilbert}).
}.

After completing his doctoral thesis, for more than 10 years served as an 
assistant of Einstein, working with him in (unsuccessful) attempts to build a 
unified field theory. 

As early as in 1917, Einstein asked Paul Ehrenfest for help to find a place for
Grommer, a ``true Russian''\footnote{
Much later a ``true Russian'' Grommer was allowed to lecture in the Belorussian
State University in Yiddish.
} 
in his words, in Russia. Later Einstein facilitated
Grommer's contacts with Russian physicists to arrange his appointment as 
professor at the Belarusian State University in Minsk in 1929. Shortly 
afterwards, Grommer was elected as a member of the Belarusian Academy of 
Sciences.

He died peacefully in 1933. Starting from 1937, he was a non-person in the 
annals of the Belarusian Academy, and in Russia in general.

Works in mathematical physics, complex analysis, analytic number theory.

\subsection{Chaim (Herman) M\"untz}
(References: 
\cite{muntz}, \cite[p.~189]{lorentz}, \cite[p.~254]{einst-files},
\cite[p.~47]{dmv}, \cite[pp.~185--187]{pinl}, \cite[p.~135--136]{ss-fleeing},
\cite{muntz-let}, \cite{muntz-let-2}, \cite{izvestiya}, \cite{freytag}).

Born in 1884 in {\L}odz (then Russia, currently Poland). Studied at Berlin with
Frobenius, Edmund Landau, and Schottky, doctoral dissertation in 1910 under 
Hermann Schwarz.

He was unable to habilitate due to some bureaucratic obstacles, and, as a 
result, never got a university position in Germany. He worked as a school 
teacher, translator, and editor of some obscure periodicals and encyclopedias. 
Having multiple interests outside mathematics, he also wrote many papers and 
several books mixing philosophy, politics, and Jewish questions.

In 1925 he pulled all the strings to get a post of professor at the newly 
established Hebrew University of Jerusalem, unsuccessfully (the chair went to 
Landau, who, however, has left back for Germany in one year).

Around 1928--1929 he briefly collaborated with Einstein (but, unlike the most
of the other Einstein's collaborators, they did not produce any joint work).

In 1929 he managed to emigrate to Russia and worked at the Leningrad University,
where in 1935 he was awarded a doctor of science degree (a Russian equivalent
of habilitation).

In 1930 he was entrusted to participate, as a key figure, in a politically 
important debate, mixing mathematical questions of intuitionism and logicism 
with ``bourgeois'' or ``Marxist'' character of mathematics.

In 1931, he managed to get a public statement from Einstein: the latter 
retracted his signature under the document condemning political trials in 
Russia, and praised Russia to the highest degree instead; after that 
Einstein, as a German civil servant, had some troubles with the German 
authorities.

In 1932, M\"untz (along with Kolman, see footnote on p.~\pageref{foot-kolman}) 
was one of the very few politically appointed Russian delegates at the 
International Congress of Mathematicians in Z\"urich.

His distinguished career in Russia came to an abrupt end in 1937, when he was 
deported on a short notice to Estonia. In his letters sent from Tallinn, he 
begged Einstein for help, to find him a visiting professor position in the 
``great democratic America''. According to Weyl, Einstein was unwilling to do 
that due to M\"untz's ``somewhat unbalanced personality''.

As a ``veteran'' emigrant, M\"untz was instrumental in bringing to Russia 
\textbf{Bergman}, \textbf{Cohn-Vossen}, \textbf{Plessner}, and \textbf{Walfisz},
and in a letter to Landau sent from Tallinn expressed worry about the fate of
all these persons.

M\"untz managed to settle afterwards in Sweden, but his mathematical work -- in
approximation theory, calculus of variations, projective geometry, and 
mathematical physics -- has been stopped at that point.

\section*{Year 1932}

\subsection{Abraham Plessner}
(References: \cite{plessner} and references therein, \cite{plessner-umn}, 
\cite[pp.~223--224]{pinl}).

Born in 1900 in {\L}\'od\'z. Studied at Giessen, G\"ottingen, and Berlin, 
completed doctorate at Giessen in 1922 under Friedrich Engel, and after that 
worked at Marburg. His habilitation submitted in Giessen in 1929 was rejected on
the pretext that he was formally a Soviet citizen. In 1932, he managed to move 
to Moscow (apparently, as a Soviet citizen, he was able to do so relatively 
hassle-free, though he was also helped by \textbf{M\"untz}), and joined the 
group of Luzin.

In 1936, he was one of the founders of \emph{Uspekhi Matematicheskikh Nauk}, and
he took his editorial work there seriously\footnote{
``Your Russian language cuts my ears'' 
({\rus ``Mne vash russki{\u\i} yazyk obrezaet ushi''}), complained Plessner 
while editing papers written by native Russian speakers 
(\cite[p.~514]{arnold-rokhlin}).
}.

In 1939 he was appointed as professor at the Moscow University, and 
simultaneously has a post in the Steklov mathematics institute. In 1949 he was
dismissed from both posts at the height of the campaign against 
``rootless cosmopolitans'' (a Russian euphemism for Jews at the time), and since
then till the end of his life he experienced a financial hardship. Around this 
time, his mathematical activity also has been stopped. Nevertheless, to his 60th
jubilee he was honored with a laudatory article in the same prestigious 
\emph{Uspekhi Matematicheskikh Nauk} he helped to establish 25 years ago.

Among his students was Vladimir Rokhlin. Israel Gelfand has called him ``a great
mathematician and a teacher''. Plessner has introduced in Moscow some areas of 
mathematics which have been not practiced there before (spectral theory, 
algebraic geometry), and was one of the founders of the Moscow school of 
functional analysis.

He died peacefully in 1961 in Moscow.

\section*{Year 1934}

\subsection{Stefan Bergman}
(References: 
\cite[pp.~17--18]{fletcher}, \cite{tomsk}, \cite[pp.~97,245]{ss-fleeing}, 
\cite{bergman}, \cite{lowner-kufarev}, \cite{tikhomirov}, 
\cite[pp.~174--175]{pinl}, \cite[pp.~16--17]{terror-and-ex}).

Born in 1895 in Cz\c{e}stochowa (then Russia, currently Poland). Got an 
engineering degree from Vienna Polytechnic in 1920, and a doctorate in 
mathematics under von Mises at Berlin in 1922. Habilitated in Berlin in 1932, 
then worked as privatdozent there. Was dismissed from Berlin by the non-Aryan 
law\footnote{
The ``Law for the Restoration of the Professional Civil Service'', introduced in
Germany on April 7, 1933. Together with another infamous ``The Reich Citizenship
Law'' introduced on September 15, 1935, made employment and normal life of Jews 
in Germany virtually impossible. These laws are referred customarily as the
``non-Aryan laws''.
} 
in 1933. He tried to get a support from the British \emph{Academic Assistance 
Council}. Despite being supported by Hadamard, he was rejected on the pretext 
that formally he was a Polish and not a German citizen.

Instead, he managed to come to Tomsk in 1934. There, in 1935, together with the
fellow emigrant \textbf{Fritz Noether} and with Theodor Molien\footnote{
Theodor Molien (1861--1941): a noted algebraist, educated in Dorpat and Leipzig,
settled in Tomsk in 1900.
}, 
he established \emph{Mitteilungen des Froschungsinstituts f\"ur Mathematik und 
Mechanik and der Kujbyschew-Universit\"at Tomsk}. In a short time, the journal
managed to attract such authors as Sergei Bernstein, Erd\"os, Khinchin, 
Kolmogorov, Kravchuk, von Neumann, and Sierpi\'nski. A famous and controversial 
paper by Einstein and \textbf{Rosen} about gravitational waves 
\cite{einst-rosen}\footnote{
Initially submitted, like a few previous papers by Einstein (and \textbf{Rosen})
of that time, to \emph{Physical Review}, a flagship American physical journal. 
As was customary for American (but less so for German) journals of that time, 
the manuscript was sent to referees, but upon receiving a (just and critical) 
referee's report, Einstein was furious, for, according to him, the manuscript 
``was sent for publication and not ... to be shown to specialists before it is 
printed''. Einstein withdrew the manuscript (``Mister \textbf{Rosen}, who was 
left for the Soviet Union, has authorized me to represent him in this matter''),
and ceased any further collaboration with \emph{Physical Review} 
(\cite[pp.~82--85,96--97]{kennefick}).
}
was republished there in 1938. Indeed, the whole enterprise looked more like a 
(first-rate) German mathematical journal on the Russian soil, and the reaction 
of authorities followed quickly: just one year after its foundation, in 1936, 
the journal and Bergman himself were under a heavy political attack during one 
of the outbreaks of the ``Luzin affair'' on the periphery\footnote{
The ``Luzin affair'' was a complex and highly controversial political campaign 
against Luzin, one of the founders of the Moscow mathematical school, charging 
him, among other things, with publishing abroad and/or in foreign languages. 
Apparently this was an (largely, unsuccessful) attempt to subordinate 
mathematics to political and ideological control (the same way as it was done, 
for example, with all humanities and with biology). At the same time, some mathematicians of the 
younger generation have tried to seize the opportunity to overthrow the ``old
guard'' and to shift the balance of powers in the Moscow (and, by extension, the
whole Russian) mathematical world. Recently a lot has been written about it, 
see, for example, \cite[\S 6]{lorentz} and \cite{neretin} (which treat the story
from entirely different viewpoints), and (numerous) references therein. Among 
the persons mentioned in this article, Aleksandrov, Khinchin, Khvorostin, 
Kolman, Kolmogorov, Schmidt, Sobolev, and Yanovskaya took part in this campaign
on the accuser's side, while Bernstein and Vinogradov on the defender's one. 
According to Aleksandrov, Luzin ``drank to the bottom of the bitter cup of vengeance of which Goethe speaks''. The whole affair has elements of ``Greek tragedy'' (\cite{simon}) or ``Shakespeare drama'' 
(\cite[p.~1]{neretin}).
}: 
the charge, unsurprisingly, was that the journal publishes in German.

In 1936, the same year, Bergman managed to escape Tomsk for Tiflis (currently
Tbilisi, Georgia), apparently rightly feeling that the situation in Tomsk is 
getting ``too hot''\footnote{
For those who managed to understand the peculiarities of life in Russia at
that times, to change the place of living, even without much hiding, was a 
common tactic to escape arrest: NKVD had quotas for how many people they should 
arrest in a given region during that or another campaign, and if they failed to 
arrest a person on their list from the first attempt, in most of the cases they
did not bother to try to reach him across Russia: it was much easier to fulfill
the quota by arresting somebody else at the same place.
}. 
He unsuccessfully advised his friend \textbf{Noether} to leave Tomsk too.

In 1937 he was forced to leave Russia, and managed to go to Paris, and then in
1939 to US, where he had a long and successful career. It seems, however, that
at the beginning in US he was not as successful in the university milieu as he
was in Russia.

Bergman was a very active person, and even during his short stays in Tomsk and 
Tiflis he managed to attract and to educate a number students: in particular, in
Tomsk he influenced Boris Abramovich Fuks\footnote{
Boris Abramovich Fuks (1907--1985): a noted specialist in complex analysis. 
``A clever person and a skillful organizer'', according to Sergei Petrovich 
Novikov.
}, and in Tiflis, Vekua\footnote{
Ilya Nestorovich Vekua (1907--1977): a distinguished Russian-Georgian 
mathematician, worked in differential and integral equations.
}. 
While in Tomsk, he arranged a visit by Hadamard in the same eventful year 
(1936). 

Works in complex analysis and differential equations.

He was said to speak all the languages of the places he lived, even for a short
period of time (including Russian and Georgian).

\subsection{Stefan Cohn-Vossen}
(References: 
\cite[p.~435]{turk}, \cite{lowner}, \cite[p.~32]{khriplovich}, 
\cite{cohn-v-umn}).

Born in Breslau (then Germany, currently Wroc{\l}aw, Poland) in 1902. 

Co-authored with Hilbert a hugely popular book \emph{Anschauliche Geometrie},
\cite{nagl-geom}.

Doctoral dissertation in 1924 at Breslau under Adolf Kneser, habilitation in 
1929 at G\"ottingen. Became a privatdozent at the University of Cologne in 1930,
and was dismissed by Germans in 1934 on the ground of the non-Aryan law.

In 1933, L\"owner recommended him for a position in US in a letter to the head 
of Dartmouth College. \textbf{Courant} recommended him for a chair in the newly 
established and ambitious Istanbul University.

He worked briefly as a gymnasium teacher in Switzerland, and emigrated to Russia
in 1934, with the help of \textbf{M\"untz} and Fritz Houtermans, an emigrant 
physicist in Kharkov. In 1935 he was awarded a doctor of science degree and was 
appointed as professor at the Leningrad University. Shortly afterwards, in 1936,
he died in Moscow from a natural cause. The official Russian obituary in 
\emph{Uspeki Matematicheskikh Nauk}\footnote{
It is, perhaps, curious to note that this was the first item in the first issue 
(under editorship of \textbf{Plessner}, among others), of what later became, 
arguably, the most influential Russian mathematical journal.
}
praised him as ``one of the most distinguished contemporary geometers''.

Works in differential geometry.

\subsection{Fritz Noether}
(References: \cite[pp.60--62]{segal}, \cite[pp.281--293]{berkovich},
\cite[p.~97]{ss-fleeing}, \cite{nkvd}, \cite{einst}, \cite{einst-1}, 
\cite{noeth-letter}, \cite{zbl}, \cite{f-noether}, \cite[\S 3]{lowner-kufarev},
\cite[pp.~203--205]{pinl}, \cite[p.~47]{terror-and-ex}, 
\cite[pp.32--33]{khriplovich}, \cite[p.~318]{roquette}.
Also, information about him is scattered over the book \cite{noether} dedicated
to his more famous sister ).

Born in 1884 in Erlangen. Brother of \textbf{Emmy Noether}. Studied at Erlangen
and M\"unich, obtained a doctor degree in 1909. Habilitation in 1911 at
Karlsruhe Technische Hochschule. During WWI, served in the German army.

Since 1922 he worked as professor at the Technical University in Breslau, from
where he was dismissed in 1934. The reason for his dismissal was twofold: on
the ground the non-Aryan law, and also due to his left-wing political views and
an open anti-Nazi stance (for example, he signed the letter in defense of 
\textbf{Gumbel}). Later the same year Noether was appointed as professor at the
Tomsk University.

In 1937, he was arrested and charged as a member of ``espionage and sabotage
group'' led by the head of the Tomsk mathematical institute Vishnevskii. His 
children were deported from Russia.

A number of prominent people did not spare efforts to rescue Noether. Weyl wrote
a letter to Muskhelishvili\footnote{
Nikolai Ivanovich Muskhelishvili (1891--1976): a Russian-Georgian mathematician,
standing high in the Russian hierarchy of that times.
}, 
trying to reach through him Beria (a chief of the Russian secret police at that
time and ethnic Georgian). Einstein wrote a letter to the Russian Commissar of 
Foreign Affairs with a request to release Noether, and helped his sons to settle
in US. Afterwards, Weyl tried to support Noether's sons financially.

Noether perished in Russia, but when and where is not exactly known. 
Fritz Houtermans, a physicist and a fellow emigrant from Germany, met him in the
infamous Butyrka prison in Moscow in 1940. The Russian officialdom claims that 
he was executed in a prison in Orel in 1941 shortly after the start of the war 
between Russia and Germany\footnote{
At least in this case the Russian officialdom has some degree of credibility,
as hundreds of prisoners were executed at that time and place in the wake of
advancing German troops.
}, 
while some witness(es?) state that he was seen as late as in the end of 1941 on the streets of Moscow.

Works mainly in applied mathematics and mathematical physics. In 1920 he 
published a pioneer paper where for the first time the index of an integral 
operator was introduced, and the first version of the index theorem was proved.

\subsection{Werner Romberg}

(References: \cite[pp.~77, 125--126]{ss-fleeing} and \cite[p.~10]{num-anal}).

Born in 1909 in Berlin. During his student days in Munich, was a member of 
Socialist Workers' Party and a staunch anti-Nazi. As early as in 1932, he was 
denied a price at a student scientific competition, despite stellar performance,
for ``lacking the necessary maturity of mind''.

His thesis adviser Sommerfeld urged him to hurry to submit his doctoral 
dissertation until it will be too late; what he successfully accomplished in 
1933. Not seeing for him any prospects in Germany, Sommerfeld, being connected 
with Russian physicists, recommended Romberg for posts in Russia. In 1934 
Romberg got a position at the Physical-Technical Institute in Dnepropetrovsk 
(currently Ukraine).

In 1937 he was forced to leave Russia, and after a brief stay in Prague, managed
to escape to Norway, first to Oslo, and later to Trondheim, where he had a long
and successful career. In 1970 he accepted professorship at the University of
Heidelberg, were he built a research school. Died in 2003.

His earlier works were in mathematical physics, and later he switched to 
numerical mathematics.

\subsection{Michael Sadowsky}
(References: 
\cite[footnote on p.~38, pp.~104,128--129,134]{ss-fleeing}, 
\cite[p.~25]{fletcher}, \cite[pp.~6--7]{terror-and-ex}, \cite{sadowsky}).

Born in 1902 in Estonia, doctoral dissertation in 1927 at Berlin under 
Georg Hamel. Habilitation in 1930 at the Technical University of Berlin, after
what he was appointed as a privatdozent there. He left Berlin in 1931\footnote{
In a few sources it is claimed that in 1931 Sadowsky has left his Berlin post 
not voluntary, but was dismissed. This seems to be wrong: Sadowsky himself was 
not Jewish and not politically active. True, he was dismissed two years later 
because of his Jewish wife, but 1931 being not 1933, the non-Aryan laws were not
yet in force, and the only reason for dismissal in 1931 could be an immense left-wing political activism, like in
the case of \textbf{Gumbel}.
}
and during 1931--1933 worked at the University of Minnesota. He returned to 
Berlin in 1933, but shortly thereafter was dismissed from the faculty of the 
Technical University because his wife was Jewish. After a brief stay in Belgium
in 1934, he went to Russia, for a short time to the Leningrad University, and 
then for 3 years to Novocherkassk, a backwater place without any scientific 
activity (it is telling that he did not publish anything during his Russian 
years).

He was expelled from Russia in 1937 on a two-days notice, and briefly stayed 
afterwards in Palestine. Neither he was able to find an academic job in 
Palestine, as he himself was ethnic Russian of the Greek Orthodox faith, and 
found the Hebrew University, the only university in Palestine of that time, ``of
an extremely nationalistic trend in the Jewish-Orthodox sense''. From 1938 he was 
in US, on the faculty of Illinois Institute of Technology.

In \cite{fletcher} it is reported on an incident when Sadowsky returned unopened
correspondence from the British \emph{Academic Assistance Council}, apparently 
angry about their unwillingness to help him.

In 1953 he joined the Renssellaer Polytechnic Institute, where he had a short --
until his untimely death in 1967 -- but distinguished career. To celebrate his
memory, the Renssellaer Polytechnic Institute established ``Michael A. Sadowsky 
Lectures in Mechanics'', and a Michael Sadowsky prize for the best Master thesis
in mechanics; some of his papers from 1930s were translated from German and 
republished recently.

Works in mechanics, in particular in the theory of elasticity, and in numerical
mathematics.

\subsection{Karl (K\'aroly) Szil\'ard}
(References: \cite[p.~281]{panorama}; also, information about him is scattered
over memoirs \cite{borin}, \cite{rumer}, and \cite{eger}).

Born in 1901 in Gy\H{o}r. Brother of \textbf{Stef\'ania Szil\'ard}\footnote{
It is frequently claimed -- for example, in \cite{borin}, \cite{panorama}, and 
\cite{panorama-szilard} -- that Karl Szil\'ard is a close relative (either 
brother, or cousin) of the famous physicist Leo Szilard. This is not true: Leo 
Szilard was born as Leo Spitz, and his parents have changed the family name to 
Szil\'ard, a common Hungarian surname, afterwards, in an apparent attempt to 
appear more Hungarian than Jewish. Also, \cite{panorama-szilard}, a very brief 
biographical sketch about Karl Szil\'ard, contains a lot of other errors.
}.

Studied at Jena and G\"ottingen. After getting a doctor degree at G\"ottingen in
1927 under \textbf{Courant}, worked in industry.

Was a member of the German communist party. Emigrated to Russia in 1934, worked
in the Central Aerohydrodynamic Institute (TsAGI) in Moscow. 

Arrested in 1938, worked in ``sharashka'' (a network of secret laboratories and
institutes utilizing forced labor of imprisoned scientists and engineers, in the
framework of Gulag, the Russian system of concentration camps) on construction 
of military aircrafts. His cellmates included his friend, the famous Russian 
physicist Yury Rumer (being released earlier then Rumer, Szil\'ard helped him to smuggle 
his scientific writings out of the prison), and the famous aircraft designers 
Andrei Tupolev and Robert Bartini. According to memoirs of that times, ``Karlusha'', as he was affectionately called by
his colleagues/cellmates, was loved by virtually everybody, including prison 
guards.

He was released in 1948, and worked again in TsAGI and in 
\emph{Referativnyi Zhurnal ``Matematika''} (a Russian analog of 
\emph{Zentralblatt} and \emph{Mathematical Reviews}).

In 1960 returned to Hungary, resumed his work in mathematics, headed the 
Department of Differential Equations of the Institute of Mathematics at 
Budapest. Facilitated contacts of his former ``sharashka'' cellmates from the 
Tupolev construction bureau with the Hungarian airline company.

Died in 1980.

Works in differential equations, complex analysis, aerodynamics.

\section*{Year 1935}

\subsection{Emanuel Lasker}

World chess champion in 1894--1921. In 1924 he was the first among chess 
international grandmasters to visit the USSR, where he got a very honorable and 
enthusiastic reception, and he praised the ``young Soviet state'' in return.

Being Jewish, Lasker was forced to leave Germany in 1933. He emigrated to Russia
in 1935 by the invitation of Nikolai Krylenko, the minister of sport and a great
chess enthusiast\footnote{
Speaks Krylenko: ``We must finish once and for all with the neutrality of 
chess ... We must organize shock-brigades ({\rus udarnye brigady}) of 
chess-players, and begin the immediate realization of a Five Year Plan for 
chess'' (\cite[p.~575 of the English edition]{souvarine}).
} 
(and concurrently the Soviet Commissar of Justice responsible for political show
trials flourishing in Russia at that time).

Lasker was immediately granted the Soviet citizenship, quickly learned Russian,
and was appointed as a coach of the Russian national chess team. He also got a 
position at the Steklov mathematics institute in Moscow. Lasker wrote a few 
seminal papers in algebra at the beginning of 1900s, but by this time he was 
over sixty and ceased to do any mathematical work long time ago, so the latter
appointment seemed to be a purely political one\footnote{
The director of the Steklov mathematics institute, a distinguished 
number-theorist Ivan Matveevich Vinogradov, mentioned many years later that 
Lasker ``enthusiastically worked on one of the mathematical problems'', without
going into details. It is safe to assume that the main occupation of Lasker at 
the institute was playing chess with Vinogradov, about what the latter had vivid
memories (see interview with Vinogradov, \cite{64}). It is also remarkable that
Vinogradov, a notorious anti-semite, apparently befriended Lasker, a Jew.
}.

Fearing political climate in Russia, Lasker left in 1937 for Netherlands, and 
then for US\footnote{
In the book \cite[pp.~177--178]{zak} another version (also repeated in a few
subsequent Russian publications) is given: Lasker, together with his wife, went
temporarily to US to visit relatives, with the intention to return to Russia in
1938. During the US trip, his wife became gravely ill, was unable to travel 
further, and so they decided to stay in US forever. Like any ``politically 
sensitive'' material of such sort published in the USSR, this version should be 
taken with a big grain of salt.
}. 
His patron Krylenko was arrested and murdered in 1938.

\section*{Year 1936}

\subsection{Myron Mathisson} (Reference: \cite{mathisson}).

Born in 1897 in Warsaw. Mathisson was a maverick, always following non-orthodox
paths in his life. Working at various odd jobs, he pursued independently 
his mathematical and physical studies, and engaged in a fruitful correspondence 
with Einstein. After Einstein's intervention, he got a doctor degree from the 
Warsaw University in 1930.

Einstein went great length in trying to arrange for Mathisson a Rockefeller 
fellowship, interacting with formal Mathisson's advisor in Warsaw, with the 
Rockefeller Foundation officers, and with not less than Rockefeller himself, all
in vein.

Being Jewish, and -- even worse -- refusing to follow the accepted ways for 
building an academic career, Mathisson did not see any prospect in Poland, and
expressed wishes to emigrate to Palestine or to Russia. However, in 1932 he 
somehow managed to habilitate at the Warsaw University, and worked as a 
privatdozent there. In 1935, by invitation of Hadamard, he lectured in 
Coll\`ege de France.

By the end of 1935, Einstein arranged an one-year visiting position for 
Mathisson at the Institute for Advanced Study in Princeton, but the letter with
invitation reached him only in 1936, when he was already in Moscow. After 
spending a short time in Moscow, he was employed as professor at the Kazan 
University. While in 1936, in a letter to Einstein, he praised his working 
conditions in Kazan, in 1937, just a year later, he complained that the 
situation there is unbearable, and the same year he fled Russia, leaving behind
him all his books and belongings.

Around the same time, his candidacy was discussed for the chair of mathematical
physics at the Hebrew University of Jerusalem. In 1938--1939 he worked at the 
Jagiellonian University in Krak\'ow, in 1939 in Paris, and in 1940 in Cambridge,
UK. He died from tuberculosis in 1940, earning an obituary note in \emph{Nature}
from Dirac, \cite{nat}, and a dedicatory paper from Hadamard.

Works in mathematical physics and partial differential equations.

\subsection{Nathan Rosen}
(References: 
\cite{ukr}, \cite{rosen-ukr}, \cite[pp.~96,192--196]{kennefick}, 
\cite[pp.146--147]{illy}).

Born in 1909 in New York, PhD in physics from MIT in 1932. Collaborator of 
Einstein (Einstein-Podolsky-Rosen paradox), was Einstein's assistant at the 
Institute of Advanced Study in 1934--1936.

Einstein was very active in efforts either to find for Rosen a suitable post, or
enable him to continue his work that or another way. One, unsuccessful, scheme 
was to do a consultant service for the Radio Corporation of America about some 
technical problems Einstein was keen about, and to use the earned fees to 
support Rosen.

Another natural possibility, owing Rosen's socialist political views, was to 
seek for him employment in Russia. Einstein asked Molotov\footnote{
Vyacheslav Molotov: a high ranking Russian diplomat, of the Molotov--Ribbentrop
pact fame.} 
for help, and this time he was more successful: in 1936, Rosen was invited to 
Kiev, and started to work at the Kiev University and the Institute of Physics of
the Ukrainian Academy of Sciences. As it often happened in such cases, his 
initial letters to Einstein were full of praise for his working conditions and 
life in Russia. However, in 1938 he returned to US, and later moved to Israel, 
where he took a number of high administrative posts in Israeli academia. 

Works in mathematical physics.

\subsection{Arnold Walfisz} 
(Reference: \cite{dmv} and references therein, \cite[p.~73]{demidovich}, 
\cite{umn}).

Born in 1892 in Warsaw. Studied in M\"unchen, Berlin, Heidelberg, and 
G\"ottingen. Defended doctoral dissertation in 1922 under Edmund Landau.
As a Polish Jew, he has no prospect of academic employment neither in Germany,
nor in Poland. Both Landau and Hardy tried, in vein, to get for him funds from 
the Rockefeller Foundation\footnote{
A negative report from the Rockefeller Foundation characterizes Walfisz as being
from ``scientifically backward country'' (\cite[p.~86]{rockefeller}); that was 
said about mathematics in the inter-war Poland. The Foundation supported 
enormously European science in the considered period, but the wisdom of its 
officers (some of them being accomplished scientists themselves) had its 
limitation.
}. 

For a while, he worked as a privatdozent at the Warsaw University, and earned 
his living in an insurance company. In 1935, together with Salomon Lubelski, he
established \emph{Acta Arithmetica}, the third specialized mathematical journal
in the world\footnote{
The first two were also Polish journals \emph{Fundamenta Mathematica} 
(set theory and logic, established in Warsaw in 1920), and 
\emph{Studia Mathematica} (functional analysis, established in Lvov in 1929). 
Prior to that, all mathematical journals were generalist ones, so the idea was 
novel and was considered with skepticism by many. Poles were forerunners in this
respect.
}.

In 1936, with the help of \textbf{M\"untz}, he managed to establish contacts 
with Georgian mathematicians, and accepted an offer from University in Tbilisi
(Tiflis at that time). He was very satisfied there -- at least at the beginning
-- with everything: working conditions, climate, colleagues. When at certain 
point a prospect arose of his possible return to Poland, due to bureaucratic 
difficulties related to his status of a foreigner in Russia, he deemed such 
outcome as a ``catastrophe''. Starting from that period, most of his works were
published in Russian in local Georgian journals. 

Walfisz established a Georgian school of analytic number theory.

He refused to speak with anybody on any ``political'' topics, and managed to 
survive unharmed the terror times in Russia. He died peacefully in 1962 in 
Tbilisi, standing high in the official Russian mathematical hierarchy --
high enough to deserve an obituary in the prestigious 
\emph{Uspekhi Matematicheskikh Nauk}.

\section*{Year 1939}

\subsection{Alfred Lustig} 
(References: \cite{lustig} and \cite[Vol. 2, p.~431]{sssr40}).

Born in 1908 in Vienna. Doctoral dissertation in physics in 1932 from the 
University of Vienna, where afterwards he worked as an assistant, and authored 
a few papers in experimental physics. After the Anschluss of Austria in 1938, 
Lustig, being Jewish, had lost his post at the university, and in 1939 was 
deported to the German-occupied part of Poland. There he managed to cross the 
border with the Russian-occupied part, and, moving further eastwards, after a 
series of odd jobs, managed to get a post at the teachers' college (later 
Pedagogical Institute) in Yelabuga in Tatarstan. There he worked till the end of
his life, with an interruption for fighting in the Russian army during WWII. 
While fighting in the army, he was heavily wounded and decorated with military orders.

Reportedly, he was found by the college authorities not qualified enough to 
teach such important subject as physics, and taught instead first German and 
then mathematics. Eventually he was appointed as the head of the department of 
mathematics. He published one paper in 1956 in mathematical analysis in a
local journal. 

Died in 1985.

Lustig spoke many languages, including the local Tatar, and was loved and 
admired by students and colleagues.

\section*{Year 1941}

\subsection{Leon Chwistek} 
(References: 
\cite[p.~97]{lapko-lyusternik} and \cite{kolm-let-chw}; also, information
about Chwistek, including a few photos, is scattered over Steinhaus' memoirs 
\cite{steinhaus}).

Born in 1884 in Krak\'ow. Served in the Polish Legions during 1914--1916. 
Painter, writer, philosopher (``A very remarkable man'', according to Mark 
Kac).

Got a doctor degree in 1922 at Krak\'ow, since 1930 worked as professor of logic
at Lvov. Being of Marxist persuasion, after occupation of Lvov by Russian 
forces, publicly praised Stalin. Left Lvov in 1941 with the retreating Russian
army\footnote{
``Professor Chwistek found a spot for himself on a Soviet lorry at the last 
minute'', witnessed Steinhaus (\cite[p.~278]{steinhaus}).
}.

While in Russia, corresponded with Kolmogorov (``A good specialist in 
mathematical logic'', attests Kolmogorov). In 1941--1943, taught mathematics at
the Tbilisi State University, after that lived in Moscow. In July 1944, 
delivered a talk at the session of the Moscow Mathematical Society.

Chwistek died in 1944 under mysterious circumstances (according to some
accounts, he died of the heart attack at the banquet in Kremlin in the presence
of Stalin; according to another ones, he was poisoned by NKVD in his residence
near Moscow).

\subsection{Zalman Skopets} (Reference: \cite{skopets}).

Born in Latvia in 1917, graduated from the university in Riga, with the 
languages of instruction being German and French. In 1941, at the outbreak of
the war between Russia and Germany, fled the advancing German army to Russia 
(Latvia was annexed by Russia in 1940, so, presumably by that time he was a 
Soviet citizen by the fact of annexation).

First worked as a school teacher, then managed to convince Russian authorities
that his Riga high education diploma is equivalent to a Russian one, and got 
accepted to graduate studies at the Moscow University. He got his candidate 
degree (a Russian equivalent of PhD) in 1946, and afterwards worked till his 
death in 1984 at the Yaroslavl Pedagogical Institute. He used to come regularly
to his native Riga, where he lectured in Latvian. He supervised many candidate 
and doctor of science theses.

Works in classical -- Euclidean and non-Euclidean -- geometry, participated in 
Kolmogorov's high school education reform.

\subsection{Mark (Marko) Vishik} 
(References: \cite{demidovich} and \cite{vishik-umn}).

His fate is similar to those of \textbf{Skopets}. 

Born in 1921 in Lvov. Studied in 1939--1941 at the Lvov University with 
\textbf{Banach}, Schauder, Mazur, Saks, and Knaster. 

In the face of the advancing German army, fled eastward to Russia. After many 
misfortunes, changing many places of residence, and many odd jobs under 
extremely hard war-time conditions, in 1942--1943 he studied at the Tbilisi 
State University, where he was helped by Muskhelishvili and \textbf{Walfisz}. 
Walfisz advised Vishik to move to Moscow, what he did. In Moscow he was 
influenced by \textbf{Plessner}, and got a candidate degree in 1947. After that
he had a distinguished career, first at the Moscow Power Engineering Institute,
and then at the Moscow University. On his 60th and 75th anniversaries he was 
honored with laudatory articles in the prestigious \emph{Uspekhi 
Matematicheskikh Nauk}.

Works in differential equations and functional analysis.

\subsection{Stanis{\l}aw Krystyn Zaremba Jr.}
(Reference: \cite[\S 2.1]{krakow} and \cite[p.~245]{steinhaus}).

Born in 1903 in Krak\'ow, son of Stanis{\l}aw Zaremba Sr., a famous Polish 
mathematician. Studied in Kr\'akow and Paris, got PhD in Vilno (then Poland, 
currently Vilnius, Lithuania), habilitated in 1936 in Kr\'akow. At the outbreak
of WWII, when Germany occupied western Poland, fled back to Vilnius, then the 
capital of independent Lithuania. After Russia occupied Lithuania, and Germans 
attacked Russia, fled eastward to Stalinabad (today Dushanbe, Tadzhikistan), 
where he worked as professor at the local teachers' college.

In 1942 Zaremba joined the Russian-supported Polish Armed Forces in the East, 
and after many adventures and a long journey through Persia, Palestine, and 
Beirut, ended up in UK, where he worked as professor at the Polish University
College in London. After the war, he worked in Madison, Quebec, and the 
University of Wales. 

Works in differential equations and stochastic processes.

\section*{Emigrations that did not occur}

\setcounter{subsection}{0}

The list is arranged in alphabetical order.

\subsection{Nachman Aronszajn}
(Reference: \cite[p.~134]{ss-fleeing}).

Aronszajn got his doctorate twice: in 1930 at Warsaw under Stefan Mazurkiewicz,
and in 1935 at Sorbonne under Fr\'echet.

In 1936, in a letter to \textbf{Courant}, Pavel Aleksandrov wrote about his 
failure to bring Aronszajn to Russia. In 1930--1940 Aronszajn was in France, in 
1940--1945 in UK, then he came back to France, and since 1948 till the end of 
his life he worked in US.

Works in functional analysis and mathematical logic.

\subsection{Stefan Banach}
(References: 
\cite[p.~66]{demidovich}, \cite{banach}, \cite{kut}, \cite{banach-umn}).

There is some evidence that a number of top Russian mathematicians (including 
Kolmogorov and Sobolev) were interested in Stefan Banach, and planned to bring 
him to Moscow and install him as a member of the Soviet Academy of Sciences. 
Banach visited Moscow in 1945 as an official guest of the Academy. Banach was 
known for a very practical apolitical approach of collaboration with 
authorities, and was loyal to the Soviet regime (as emphasized in his obituary
published in \emph{Uspekhi Matematicheskikh Nauk}, during the Russian occupation
of Lvov in 1939--1941, he was installed as a dean of the physical-mathematical
faculty, and a member of the city council). On the other hand, after the 
end of the war he was offered a chair in Kr\'akow, while at the same time some 
Polish sources claim that he intended to move to Wroc{\l}aw, where most of the 
few remaining people from the Lvov University were transferred to. In any case,
Banach was mortally ill at that time, and died the same year, 1945.

\subsection{Richard Courant}
(References: 
\cite{kolm-let}, \cite[p.~219]{john}, \cite{reid-courant}).

From the 1932 letter of Kolmogorov to Pavel Aleksandrov: ``By 1936 ... according
to Kolman, Courant will sit in Moscow''\footnote{
In his memoirs Kolman writes that Courant has approached him during the 1932
International Congress of Mathematicians in Z\"urich with an inquiry about 
possibility to settle in Russia 
(\cite[pp.~237--238 of the 2011 edition]{kolman}). Courant did not consider 
emigration until his dismissal in 1933 by the non-Aryan law which took him by
surprise, and which he tried to appeal, and definitely had no ``invitation from
the New York University'' in 1932 as Kolman writes. So the whole story is not 
clear (Kolman is an unreliable source, while Kolmogorov is a reliable one), but it is reasonably to assume that
some interest on the part of Courant to Russia did exist at that time.
}.

From the memoir of Fritz John about events in 1933: ``Courant pursued various 
leads, for example to Istanbul or Odessa, which came to nothing (fortunately, as
it turned out)''.

\subsection{Werner and K\"ate Fenchel}
(Reference: \cite[footnote 95]{denmark}).

Werner Fenchel was a noted mathematician working in geometry and optimization
theory. His wife, K\"ate, n\'ee Sperling, was a mathematician working in group
theory.

Werner Fenchel was dismissed by the non-Aryan law from his post at G\"ottingen 
in 1933, and the couple fled to Denmark, were they remained till the end of 
their lives (with a brief escape to Sweden in 1943--1945 as a part of the mass 
rescue company of Danish Jews). In Denmark Werner Fenchel had a distinguished 
career, while K\"ate worked for a few years as a secretary of Harald Bohr.

However, during the 1930s their future was not certain at all, and Fenchels 
applied for support to the British \emph{Society for the Protection of Science 
and Learning}, indicating Russia as one of their possible destinations. Yet 
there is seemingly no evidence that they made any serious attempt to settle in 
Russia.

\subsection{Emil Julius Gumbel}

(References: 
\cite[pp.~117--118,202--203]{rockefeller}, \cite{sheynin}, 
\cite[pp.~62--63]{segal}, \cite[pp.~255,262]{einst-files}, 
\cite[\S 1]{gumbel-stat}).

Staunch pacifist, human rights activist, and anti-Nazi.

In 1925 Gumbel approached Einstein, asking for Russian contacts which would 
facilitate for him a possibility to get a post in Russia. He visited Russia in 
the winter of 1925--1926, worked there on ``mathematical archives of Karl 
Marx''\footnote{
``Karl Marx's mathematics'' was, for a long time, a popular speciality among a 
number of Russian philosophers, historians, and mathematicians, including Kolman
and Yanovskaya.
}, 
and widely publicized his positive views of the country. He briefly visited 
Russia again in 1932.

Gumbel outspoken views and political activism were too radical even for moderate
anti-Nazi supporters of the Weimar Republic, and his political opponents 
constantly tried to remove him from his modest university posts in Germany. In 
1931 a public letter for Gumbel's defense was signed, among others, by Einstein,
and by \textbf{Emmy} and \textbf{Fritz Noether} (but not by \textbf{Courant}, who 
refused to sign). Eventually, Gumbel lost his post at Heidelberg in 1932, and 
the same year came to Paris, by \'Emile Borel's invitation, to work at the 
Institute Henri Poincar\'e, and later was supported by CNRS. He was rejected several times by the Rockefeller Foundation, apparently to a large 
degree due to his political activism. On many occasions, Einstein tried to help
him to secure posts in different countries. Finally, due to efforts of Einstein
and some leading American statisticians, Gumbel managed to come to US in 1940.

Still, even in US, as well as in post-war Germany, he was often not welcomed and
got rejected on many occasions.

Works in statistics.

\subsection{Wladimir Lewicki (Vladimir Iosifovich Levitskii)\protect\footnote{
Spelled as W{\l}odzimierz Lewickij in \cite{duda}, and as Wl. Lewicky and 
W. Lewickyj in \emph{Zentralblatt}.}}
(References: \cite[pp.~27--28,151,184]{duda}, \cite{khobzei}, \cite{krav-let}).

Doctorate in 1901 at the Lvov University. He can be compared with 
\textbf{Czajkowski}, a fellow Ukrainian in Lvov at that times: like 
\textbf{Czajkowski}, Lewicki worked in Lvov at the underground private Ukrainian
university and in high schools; at the same time, unlike \textbf{Czajkowski}, he
maintained close contacts with the Polish Lvov school of mathematics, and 
published, at least part of his papers, in German. In 1929 Kravchuk encouraged 
Lewicki to emigrate to Ukraine, offering positions either in Kiev or in Kharkov,
but nothing came out of these efforts.

Upon Russian occupation of Lvov in 1939, and its subsequent ukrainization, 
Lewicki was appointed as professor at the Lvov University. After WWII, upon 
exodus of Polish mathematicians, he occupied the chair of \textbf{Banach}, and 
tried to preserve what little is remained from the traditions of the Lvov 
mathematical school. Died peacefully in 1956.

Works in the theory of analytic functions.

\subsection{Rudolf Karl L\"uneburg}
(Reference: \cite[pp.~134,166,180--185]{ss-fleeing})

Doctoral dissertation in 1930 at G\"ottingen. Was dismissed in 1933 from his
assistant post at G\"ottingen for political reasons.

In 1935, Pavel Aleksandrov wrote in a letter to \textbf{Courant} about slim 
chances to find a post for L\"uneburg in Russia: the reason, according to 
Aleksandrov, was that L\"uneburg publishes less than other mathematicians of his
caliber\footnote{
An evidence that the so-called bibliometrics played a role also in those earlier
years! As an exercise in an alternative history, one may speculate on what 
impact all these ``Impact factors'' and ``Hirsch indices'' would have on the 
survival of that or another mathematician in the Nazi era, would they be 
invented a little bit earlier.
}.

In 1934--1935 L\"uneburg worked at Utrecht and Leiden. In 1935 he managed to 
come to US. Due to his conflicts, going back to their G\"ottingen days, with 
such influential emigrants as Weyl and Busemann, he failed to find a job in 
academia, and worked in industry until his death in 1949 in an automobile 
accident.

Works in function theory and mathematical optics\footnote{
In \cite{ss-fleeing}, L\"uneburg is wrongly described as a topologist.
}.

\subsection{Kurt Mahler}
(References: \cite[footnote on p.~132]{ss-fleeing} and 
\cite[p.~165]{nossum-kotulek}).

Kurt Mahler, a distinguished number-theorist, got his doctor degree in 1927 at
Frankfurt under Siegel. Before assuming his first post as a privatdozent at 
K\"oningsberg in 1933, he decided to leave Germany. He managed to find some 
funding in UK, but his situation was not secure there, and around 1935 he had a
hope to get a position at the Saratov University, with the help of 
Khinchin\footnote{
During that time the Saratov University flourished under the rector Gavriil 
Kirillovich Khvorostin, an alumni of the Moscow University, visitor to Berlin 
and G\"ottingen, and an active and ambitious administrator with connections at 
the top of Soviet bureaucracy. During his short reign in 1935--1937, Khvorostin 
managed to attract to the mathematics department of the Saratov University such
first-rate people as Israel Isaakovich Gordon, Aleksandr Yakovlevich Khinchin, 
Aleksandr Gennadievich Kurosh, Ivan Georgievich Petrovskii, and Viktor 
Vladimirovich Wagner (about whom and similar people he once quoted verses of the
Russian poet Sergei Esenin: ``It is hard to handle the cattle with a twig''; the
Russian original ``{\rus Trudno khvorostino{\u\i} upravlyat{\cprime} 
skotino{\u\i}}'' is based on a wordplay: Khvorostin vs. ``khvorostina'', a 
twig). The glorious period of ``G\"ottingen on Volga'', as Khvorostin liked to 
put it, came to end in 1937, when Khvorostin was arrested and later murdered, and most of the scientists hired by him were either fired, 
or left Saratov voluntary fearing persecutions (\cite[pp.~287--288]{umn-chronicle}, \cite[pp.~19--21]{gordon}).
}.
However, nothing came out of these efforts, and Mahler established himself in UK
and later in Australia.

\subsection{Emmy Noether}
(References: \cite{noether}, including \cite{alexandrov}, \cite{roquette},
\cite[pp.~67--72]{transcending}).

Emmy Noether has spent the 1928--1929 academic year in Moscow, were she 
influenced Pavel Aleksandrov, Pontryagin, Schmidt, and others. She enjoyed this
visit very much.

After the introduction of the non-Aryan laws in Germany, Noether's friend 
Aleksandrov tried to arrange for her a chair in Moscow, unsuccessfully. 
(Earlier Noether helped Aleksandrov to obtain a Rockefeller fellowship, and 
recommended him for a professorship at Halle). Instead, in 1933 she departed to
US. After her untimely death in 1935, she was honored in a memorial session of 
the Moscow Mathematical Society, attended by her brother \textbf{Fritz}.

\subsection{Felix Pollaczek} 
(References: 
\cite{pollaczek} and references therein, \cite[pp.~111,133]{ss-fleeing},
\cite[p.~122]{rider}).

Doctoral dissertation in 1922 at Berlin under Issai Schur. Pollaczek was a 
versatile mathematician and electrical engineer, working in number theory, 
mathematical analysis, probability, mathematical physics, and other areas. In 
1933 he was dismissed by the non-Aryan law from his post at the German Postal, 
Telephone and Telegraph Service. He moved to Paris, but was forced to leave 
France in 1936 as his visa was not prolonged. Around this time Khinchin and 
omnipresent \textbf{M\"untz} have tried to bring Pollaczek to Tiflis or Baku 
(currently Azerbaidzhan). These efforts failed\footnote{
In the first issue of \emph{Uspekhi Matematicheskikh Nauk} it is reported that
the department of algebra in the mathematical institute at Tiflis is headed
by ``Professor Pollaczek from France'' (\cite[p.~287]{umn-chronicle}), on what
Kolmogorov, later criticizing the quality of information in the newly 
established journal, sarcastically remarked that ``the information preempted the
event'' (\cite{umn-kolm}).
};
Pollaczek was able to come to Russia, by invitation of Khinchin, only for a 
short visit in 1937.

Instead, in 1938 Pollaczek managed to move back to France. Such people as 
von K\'amr\'an, Veblen, and Weyl did not spare efforts to bring him to US, all 
in vein. In 1942 Pollaczek managed to get an appointment in Lima, but was unable
to move out of the chaotic war-time France. He survived during the German 
occupation, and remained in France till the end of his life. After the war his
only source of income was a meager CNRS stipend, so he experienced a constant
financial hardship, but continued to work fruitfully nevertheless, writing many
more papers.

\subsection{Erich Rothe and Hildegard Rothe-Ille}
(References: \cite[p.~133]{ss-fleeing}, \cite[pp.~208--209]{pinl}, 
\cite[pp.~13--14]{ille}).

Erich Rothe has defended a doctoral dissertation in 1927 at Berlin under Erhard
Schmidt and Richard von Mises. He habilitated in 1928 at the Technical 
University of Breslau, where he worked as an assistant of 
\textbf{Fritz Noether}; in 1931--1935 he worked at the University of Breslau,
from where he was dismissed in 1935 by the non-Aryan laws.

In 1934, Oswald Veblen in a letter to Richard Brauer advised Rothe to seek 
employment in Russia, where, according to Veblen, he has much more chances than
in US. Eventually, Rothe went to US in 1937, where, after a number of temporary 
posts, he had a long and successful career at the University of Michigan.

Works in partial differential equations and functional analysis.

His wife Hildegard Rothe-Ille has defended a doctoral dissertation in 1924 at
Berlin under Issai Schur. She followed her husband to US, where she taught at a
small private college, and died in 1942 from cancer. Besides her doctoral  
dissertation, she published only one short paper in 1926 about 
positive-definite polynomials.

\subsection{Hans Schwerdtfeger}
(References: \cite{aczel}, \cite{born}, \cite{einst-born}, 
\cite[pp.~123,137,139,165]{ss-fleeing}).

Doctoral dissertation in 1934 at Bonn under Toeplitz.

Schwerdtfeger's friend \textbf{Cohn-Vossen} tried to secure for him a place in 
Russia. These efforts were brought to halt due to \textbf{Cohn-Vossen}'s death. 
Max Born also tried to use his contacts with Russian physicists to find him a 
place in Russia, and urged Einstein to write on behalf of Schwerdtfeger to 
Molotov or Schmidt. Einstein refused: he was of quite a low opinion of 
Schwerdtfeger as a scientist, and would not like to compromise his reputation by
recommending a ``mediocrity''. Weyl was also of a low opinion of Schwerdtfeger,
and refused to help him to settle in US.

Schwerdtfeger was a staunch anti-Nazi, and in 1936 he was forced to flee the 
country to Prague (by Born's version, he went to Prague to facilitate contacts 
with \textbf{Cohn-Vossen}), then to Switzerland (facing there a threat to be 
deported back to Germany). After a long journey involving several other 
countries, in 1940 he settled, with the help of Max Born, in Australia, and 
later in Canada.

Works in algebra and complex analysis.

\subsection{Wolfgang Sternberg}
(References: 
\cite{sternberg-obi}, \cite[pp.~128--129,134,209]{ss-fleeing}, 
\cite[pp.~216--217 of the 1996 edition]{reid-courant}, 
\cite[pp.~209--211]{pinl}).

Doctoral dissertation in 1912 at Breslau, habilitation in 1920 at Heidelberg.

In 1935 Sternberg was dismissed by the non-Aryan laws from his post at Breslau,
and briefly went to Palestine. However, in Palestine he was unhappy, due to 
various factors, such as ``Jewish nationalism'', lack of knowledge of Hebrew, 
and bad climate. He worked at the Hebrew University without renumeration, and 
actively sought possibilities to emigrate to other places, among them to Russia.

In 1935 he reported to his classmate \textbf{Courant} about a failed attempt to
get a post in Russia, through the Leningrad mathematician Vladimir Ivanovich 
Smirnov, and biologist and an emigrant from Germany Julius Schaxel.

In the same 1936 letter mentioned under \textbf{Aronszajn}, Aleksandrov wrote to
\textbf{Courant} that he is dubious about Sternberg's prospects in Russia, as 
earlier he was unable to do anything for such ``outstanding'', in his words, 
people as \textbf{L\"uneburg}, \textbf{Zorn}, and \textbf{Aronszajn}.

Sternberg left Palestine for Prague and eventually, in 1939, managed to reach
US where he had a great difficulty to find an employment, despite help from 
\textbf{Courant}. Eventually, he went through a number of temporary jobs at the
Cornell University and other places, until his retirement in 1948.

Works in potential theory and integral equations.

\subsection{Dirk and Saly Ruth Struik} 
(References: \cite{rowe}, \cite{alberts}, and 
\cite[p.~188--189 of the 2011 edition]{kolman}).

Dirk Struik was a staunch Marxist and a member of Dutch communist party (at a
certain point in his life, being criticized by his fellow party members for 
studying ``all that bourgeois science'', he deliberated whether to continue his
mathematical studies, or to become a party functionary). Due to his political 
views, he had difficulty to get a permanent academic appointment anywhere in 
Europe. 

His brother Anton, also a communist, went to Russia as early as 1921, where he 
worked as an engineer until return to Netherlands in 1930. This experience, 
quite probably, motivated Dirk to follow the steps of his brother and to settle
in Russia.

In 1926, he had to make a difficult choice between two offers: from Otto 
Yulievich Schmidt in Moscow, and from Norbert Wiener in MIT (whom he befriended
during his stay in G\"ottingen). After much deliberation, he had choose the 
latter, having in mind to accept the Russian offer sometime in the 
future\footnote{
In the 1987 interview Struik admitted that, would he be settled in Russia, his 
``natural Dutch obstinacy ... might have gotten in the way and brought [him] 
into conflict'' with Russian authorities (\cite[p.~24]{rowe-interv}).
}.
He remained at MIT till the end of his long career (where he had troubles during
the McCarthy era). He visited Russia briefly only in 1934, where he participated
in the famous international Moscow conference on tensor and differential 
geometry, and stayed with his friend Kolman\footnote{
According to Kolman's memoirs, which should be taken with a big grain of salt.
}.

Works in differential geometry and history of mathematics.

Dirk Struik's wife, Saly Ruth Struik (n\'ee Ramler), also a mathematician, got 
a doctorate at Prague in 1919 under Gerhard Kowalewski and Georg Pick. Besides 
her doctoral dissertation, she authored a single paper on affine geometry, wrote
commentaries to the Italian edition of Euclid's ``Elements'', and coauthored a 
couple of papers with her husband about history of mathematics.

\subsection{Ludwig Wittgenstein}
(References: Wittgenstein's letters 
\cite[letters 190--193 to/from John Maynard Keynes, July 1935]{wittgenstein}; 
also \cite{moran}, \cite{omahony}, \cite[Chapters 16 and 17]{monk}, 
\cite{biryukov}, and memoir of Fania Pascal, his teacher of Russian 
\cite[pp.~30--31, 37--38 of the original 1973 publication]{pascal}. These 
accounts are often contradictory, and accuse each other in wrong interpretation
of Wittgenstein's intentions, ideas, and works).

Wittgenstein interest in Soviet Russia started, probably, as early as in the 
beginning of 1920s, and intensified in the beginning of 1930s due to influence 
of his friends Russell and Keynes\footnote{
Russell was overly critical about Soviet Russia; as suggested in 
\cite[p.~248]{monk}, ``if Russel hated it so much there must be something good
about it''. Keynes presented Russian Marxism as a religious faith which, again,
according to \cite{monk} could attract Wittgenstein.
}.
During 1933--1935 Wittgenstein learned the Russian language, and tried to secure
the necessary contacts -- including multiple contacts with communists in 
Cambridge -- which would allow him to make a trip to Russia, and, ultimately, to
get a permanent residence and employment there.

Apparently Wittgenstein's desire to go to Russia was so great that he, in the
normal circumstances not hesitating a moment to use harsh words towards his 
colleagues and friends, and behaving generally in an extremely eccentric way, 
demonstrated a lot of uncharacteristic for him humility: in a letter to his 
friend with Russian connections he begged to convince Russian authorities
that he is ``in no way politically dangerous'', and during the interview with 
the Russian ambassador in London wore a suit, for the first (and, perhaps, the 
last) time in many years.

He traveled in 1935 to Leningrad and Moscow for 3 weeks in an attempt to 
secure employment in Russia, and had encounters with Sofia Yanovskaya\footnote{
Sofia Aleksandrovna Yanovskaya (1896--1966; spelled as Janovskaja or Janovskaya
by the cited English-speaking Wittgenstein scholars), a well-known and 
controversial figure in the Moscow mathematical milieu between 1930s and 1960s.
Professor of logic and history of mathematics at the Moscow University, a 
combatant in the Russian civil war, and a watchdog of ``Marxist-Leninist 
character of logic and mathematics''; at the same time she was credited for 
saving mathematical logic in Russia from even more zealous watchdogs, and was
praised as a good saint helping a few brilliant young mathematicians to 
establish themselves in Moscow. Criticized Schmidt for ``idealism'' and 
deviating from ``communist party analysis'' in his works in group theory, and 
advised Wittgenstein ``to read more Hegel''.
}.
According to some reports, Wittgenstein did not want to be employed there in 
academia, but had vague ideas ``to practice medicine in Russia'', or to be 
somehow involved in ``the newly colonized parts at the periphery of the 
USSR''\footnote{
This is, perhaps, not surprising, as at the different moments of his hectic 
life, Wittgenstein worked as a soldier, as a gardener's assistant, as a 
packing-cases maker, as an elementary school teacher, and as a hospital porter.}. 
Instead, he was offered an employment at the philosophical department either in 
Moscow or in Kazan, but shortly thereafter the Russian authorities have changed
their mind and retracted the offer. 

Initially Wittgenstein planned to go to Russia together with his homosexual 
partner Francis Skinner\footnote{
A promising student of mathematics at Cambridge in early 1930s. Under 
Wittgenstein's influence, switched from mathematics to philosophy, and then
abandoned academia altogether. Skinner is not qualified as mathematician 
according to our standards, and therefore does not deserve a separate entry 
under ``Emigrations that did not occur''.
}, 
but the plans failed due the serious illness of the latter. This, according to 
some sources, could be an additional reason why Wittgenstein has abandoned his 
plans to settle in Russia.

\subsection{Max Zorn} (Reference: \cite[p.~134]{ss-fleeing}).

A famous mathematician (of the Zorn lemma) working in algebra and numerical 
analysis, doctorate in 1930 at Hamburg under Emil Artin. Zorn worked till 1933 
at the University of Halle. In 1933, he was denied habilitation and was forced 
to leave Germany, being actively opposed to the Nazi regime.

In the same 1936 letter to \textbf{Courant} mentioned under \textbf{Aronszajn}
and \textbf{Sternberg}, Aleksandrov wrote that he was unable to bring Zorn to 
Russia.

In 1934 Zorn emigrated to USA, where he had a long and distinguished career.

\section*{Conclusion}

Why there were so few emigrations? One of the reasons was a total absence in 
Russia of any central administrative body supervising the situation and helping
the dismissed scholars to find financial sources and jobs. It is somewhat ironic
that with all its self-praised centralized planning, Russia has nothing like 
American \emph{The Emergency Committee in Aid of Displaced European Scholars}, 
or British \emph{Academic Assistance Council} (renamed in 1936 to 
\emph{Society for the Protection of Science and Learning}), or French 
\emph{Comit\'e des Savants} (operational till the occupation of France by 
Germany in 1940), or Swiss \emph{Notgemeinschaft Deutscher Wissentschaftler 
im Ausland}. In Russia, unlike all these countries, each case was processed on
ad-hoc basis: interested mathematicians (like Aleksandrov), or administrators 
(like Khvorostin), or political figures (like Krylenko) were persuading their 
respective political bosses to grant a permission to bring to the country that 
or another foreign scholar.

Another, related, reason was a byzantine character of the Russian bureaucracy: 
to employ a foreigner was a ``political'' decision, requiring involvement on a 
higher and higher levels of the Soviet hierarchy (it is reported in some sources
that in some cases the decision to offer or to reject an employment in the 
country was taken by Stalin). Apparently the political and bureaucratic 
obstacles were less at the Russian periphery than in the center: among the 22 
emigrants, less than a half were able to settle, if only temporarily, at Moscow
and Leningrad, the two Russian great scientific centers, and the rest found 
employment at the periphery. One may speculate that if Aleksandrov would lobby
the authorities for \textbf{Emmy Noether}'s employment not at Moscow, but at 
Tbilisi, Minsk, or Tomsk, he might be more successful. Perhaps, it is also not 
coincidental, that most of those very few emigrants who managed to leave a trace
in Russian mathematics, at the time of their emigration had, in that or another 
way, certain rights to the Soviet citizenship: either being born on what was 
Russian territories in the past, or became Soviet citizens by the fact of 
Russian annexation of new territories.

One might think that another obvious reason was a troublesome, to put it mildly,
political situation in Russia. However, it appears that in most of the cases 
this was not a decisive factor in choosing a possible country for emigration; 
the true nature of the Russian political regime became apparent to people to the
full extent only in the hindsight, after they spent some time in the country. In
the considered period, end of 1920s--1930s, Russia appeared as an attractive 
destination for many people, especially educated intellectuals with left-wing 
political leaning: a country not without its rough edges, but successfully 
building a progressive new society, a place where scientists are respected and highly rewarded for their research and pedagogical 
endeavors\footnote{
Evidences of this abound, here is just a couple of quotes from the already cited
sources: ``... all free minded people look at Russia with certain admiration and
with much interest'', wrote Richard Brauer to Szeg\"o in 1934 
(\cite[p.~133]{ss-fleeing}). ``Go back to Lw\'ow. The Bolsheviks idolize 
professors. They won't harm you'', advised to Hugo Steinhaus an officer in 
September 1939 on the Polish-Hungarian border, when Steinhaus was deliberating 
whether or not to flee the advancing Russian army (\cite[p.~224]{steinhaus}).
}. 
Also, at least for a certain period of time, with all the post-WWI devastation 
and later Great Depression, the working and economic conditions of scientists in
Russia seemed to be at least not worse then those of their colleagues in Europe\footnote{
In 1927 Sergei Bernstein, in a talk with an officer from the Rockefeller 
Foundation, claimed that, speaking of Russia, ``the professor's situation from 
the material point of view was comparable to that of the French professors; 
mathematicians are naturally given all freedom in their work'' 
(\cite[p.~121]{rockefeller}). Bernstein, on several occasions, courageously 
voiced his opinion against Russian authorities, and, by all accounts, this 
should be taken as a sincere description of the situation and not as a communist
propaganda.
}.
The British \emph{Society for the Protection of Science and Learning} in 
desperate efforts to help the growing number of displaced scholars, encouraged 
them to apply for jobs at Leningrad, among other ``exotic'' places 
(\cite[p.~139]{rider}).

With the exception of \textbf{Stef\'ania Bauer}, the first wave of emigration 
started in 1929, despite that some people (\textbf{Czajkowski}, 
\textbf{Grommer}, \textbf{Gumbel}) were seeking opportunities to emigrate to 
Russia earlier. This is not accidental: 1929 is the year of consolidation of 
Stalin's power, and the start of the radical changes in the economic policy -- 
the so-called ``Great Turn'' -- a drastic acceleration of forced 
industrialization of the country, and adoption of five-year plans of economic 
development. This, in its turn, required a significant improvement of education,
with focus on technical and mathematical education, a task hindered by lack of qualified mathematics teachers in universities (see, for example, a 
talk by Otto Yulievich Schmidt at the first USSR mathematical symposium in 1930,
\cite[p.~58]{lapko-lyusternik}). It is reasonable to assume that, from the point
of view of Russian authorities of that time, the emigration of scientists could 
help in all these noble goals. 

The next batch of emigrants came in 1934--1936, after the wave of dismissals in 
Germany which started in 1933. The emigration has stopped entirely in 1937, the
year of ``Great Purge'', and has resumed, if only slightly, at the beginning of
WWII (invasion of Poland by Germany and Russia in 1939, and invasion of Russia 
by Germany in 1941), under entirely different circumstances.

One characteristic feature is a frequent appearance of Einstein as a facilitator
between potential emigrants and Russian authorities\footnote{
In some recent Russian mass media and popular literature, one can find claims 
that Einstein himself was offered and even considered posts in Russia. One of 
the most popular repeated accounts goes as follows: in 1935, Khvorostin has 
offered to Einstein a position at the Saratov University. Einstein refused on 
the pretext that he never will be able to master Russian. Khvorostin countered 
with a fantastic plan to establish an Academy of Volga Germans, with German as 
an official language, and Einstein, stationed in Saratov, as the head of the 
Academy. Neither Einstein archives in Jerusalem and Princeton, nor any other 
sources do seem to contain anything corroborating the story; in 1935 Einstein 
was firmly established at the Institute of Advanced Study, and did not seek 
another positions. Fantastic as it is, this story probably reflects Khvorostin's
high aspirations as mirrored in the public consciousness of later times.
}.
This is hardly surprising: Einstein was a world celebrity with a huge number of
connections, usually was sympathetic to other people's problems, and, being 
liberal with left-wing leaning, had quite a soft spot for the communist 
Russia\footnote{
As early as in 1923, Einstein was one of the founders of the ``Association of 
Friends of the New Russia''. Also, see an incident with his public statement 
described in connection with \textbf{M\"untz}, and footnote on 
p.~\pageref{foot-eins}.
}.
Another important facilitator on the Russian side was Pavel Aleksandrov. This is
natural either: as a frequent visitor to the European mathematical centers in 
1920s, Aleksandrov had many personal friends in the West; on the other hand, he
was a patriot of his country, and apparently genuinely wanted to strengthen 
Russian mathematics by bringing excellent people from abroad.

Russia had a long tradition, starting from the St. Petersburg Academy of 
Sciences created by Peter the Great, to invite foreign scholars, mathematicians
in particular, and benefit greatly from their presence in the country. This 
time, however, the political climate and overall situation were different. A 
crude ``statistics'': out of 22 mathematicians who emigrated to Russia between 
1925 and 1941, three were murdered (\textbf{Bauer}, \textbf{Burstin}, 
\textbf{Noether}), one died mysteriously (\textbf{Chwistek}), four were deported
or forced to leave the country (\textbf{M\"untz}, \textbf{Bergman}, 
\textbf{Romberg}, \textbf{Sadowsky}), four have left themselves after a short 
period of time (\textbf{Lasker}, \textbf{Mathisson}, \textbf{Rosen}, 
\textbf{Zaremba}), four were prosecuted in that or another way, being 
imprisoned, or deprived from academic employment, or deported from the center to
periphery (\textbf{Czajkowski}, \textbf{Frankl}, \textbf{Plessner}, 
\textbf{Szil\'ard}), and in one case all the memory about a person was, after 
his natural death, erased from the official history (\textbf{Grommer})\footnote{
Another favorite action of NKVD, especially accelerating during 1939--1941, the 
time of Russian-German ``non-ag\-res\-si\-on treaty'', was transferring German 
refugees to NKVD's colleagues from Gestapo; for example, this is what happened 
to Fritz Houtermans, a physicist mentioned above in connection with 
\textbf{Cohn-Vossen} and \textbf{Fritz Noether}. This fate, however, was not 
shared by anyone of emigrant mathematicians. (Despite doing some, apparently 
serious -- it met approval of van der Waerden -- number theory in his head while
in NKVD captivity, Houtermans still does not qualify for being mathematician, 
and is not considered by us here separately).
}. 
It seems that among those who emigrated before 1937, only \textbf{Frankl}, 
\textbf{Plessner}, and \textbf{Walfisz}, and, to a lesser degree (but quite 
amazingly, taking into account a very short period of his activity in the 
country), \textbf{Bergman}, were able to influence the local mathematical 
community in some substantial way. In the last small group of emigrants who came
after 1939, \textbf{Skopets} and \textbf{Vishik} also have benefited 
significantly their new country of residence.

On the other hand, some emigrants (\textbf{Burstin}, \textbf{Frankl}, 
\textbf{M\"untz}) have participated, if not actively, in political campaigns 
aimed to subordinate mathematics to the communist ideology, and targeting their
colleagues.

The official Russian historiography of that period followed the general trend in
the country to pretend that ``non-desirable'' persons do not exist. After the
war, two big reference works were published, presumably listing all the Soviet 
mathematicians, and all their published works for the period from 1917 till 1947
and 1957, respectively (\cite{sssr30} and \cite{sssr40}). As expected, those 
mathematicians who had left the country, or were imprisoned, or executed, are 
not mentioned there, though there are some inexplicable exceptions (thus, in 
both volumes \textbf{M\"untz} and \textbf{Romberg} are mentioned).

Of course, in the considered period Russian mathematics was already a first-rate
one and was developing successfully without ``help from abroad''. Some authors, 
describing the fate of that or another person, expressed a relief that his or
her emigration to Russia did not materialize (see, for example, a remark by 
Fritz John about \textbf{Courant}): indeed, wouldn't it be quite disheartening 
to see either \textbf{Courant}, or \textbf{Emmy Noether}, or 
\textbf{Pollaczek}, or \textbf{Struik} to suffer under the Russian dictatorship?
At the same time, looking on life and deeds of the emigrants to US, the persons
who, to a large extent, shaped the American mathematical culture: Artin, 
Ahlfors, Bargmann, Bers, Bochner, two Brauers, Busemann, Carnap, Chevalley, 
\textbf{Courant}, Eilenberg, Feller, Friedrichs,  Hurewicz, John, 
von K\'arm\'an, Lewy, Loewner, Menger, von Mises, Neugebauer, von Neumann, 
Neyman, Polya, Prager, Rademacher, Rad\'o, Schoenberg, Szeg\H{o}, Tarski, 
Taussky-Todd, Uhlenbeck, Ulam, Wald, Warschawski, Weil, Weyl, Wigner, Wintner, 
Zariski, \textbf{Zorn}, Zygmund, and the scores of others, one cannot help but 
wonder how differently two superpowers used the opportunities provided by exodus
of brilliant mathematical minds from Europe at that times, and what amazing 
possibilities for Russian mathematics were lost.

\section*{Acknowledgements}

Thanks are due to Willi J\"ager, Eugen Paal, Vadim Rudnev, 
Witold Wi\c{e}s{\l}aw, Maria Zusmanovich, and especially Galina Sinkevich, for 
useful remarks and/or help with literature; and to Albert Einstein Archives at 
the Hebrew University of Jerusalem for providing copies of (yet unpublished)
letters and other documents.

\end{document}